\documentclass[a4paper,12pt]{article}
\usepackage{hyperref}
\usepackage{amssymb}
\usepackage{amsmath}
\usepackage{amsthm}
\usepackage[latin1]{inputenc} 
\usepackage{ulem} 
\usepackage{tikz} 
\usetikzlibrary{shapes.misc}
\usepackage{imakeidx}
\makeindex
  
\newtheorem{defi}{Definition}[section]
\newtheorem{prop}[defi]{Proposition}
\newtheorem{thm}[defi]{Theorem}
\newtheorem{lem}[defi]{Lemma}
\newtheorem*{lem*}{Lemma}
\newtheorem{lemma}{Lemma}
\newtheorem{cor}[defi]{Corollary}
\def\R{{\mathbb R}}
\def\C{{\mathbb C}}

\def\N{{\mathbb N}}

\def\eps{{\varepsilon}}

\newcommand {\tx}[1] {\textrm{#1}}

\newcommand{\pair}[1]{\langle{#1}\rangle}

\newcommand {\link}[1]{\href{#1}{link}}

\begin{document}

\title{K-theory of group Banach algebras and Banach property RD}
\author{Benben Liao and Guoliang Yu}

\newcommand{\Addresses}{{
  \bigskip
  \footnotesize

  Benben Liao, 
  \textsc{Department of Mathematics, Texas A\&M University, College Station, Texas, USA,} \& \textsc{Shanghai Center for Mathematical Sciences, Shanghai, China}
  \par\nopagebreak
  \textit{E-mail address}: \texttt{liaob@math.tamu.edu}

  \medskip

  Guoliang Yu, \textsc{Department of Mathematics, Texas A\&M University, College Station, Texas, USA,} \& \textsc{Shanghai Center for Mathematical Sciences, Shanghai, China}
  \par\nopagebreak
  \textit{E-mail address}: \texttt{guoliangyu@math.tamu.edu}
}}

\maketitle

\abstract
{We investigate Banach algebras of convolution operators on the $L^p$ spaces of a locally compact group, and their K-theory. We show that for a discrete group, the corresponding K-theory groups depend continuously on $p$ in an inductive sense. Via a Banach version of property RD, we show that for a large class of groups, the K-theory groups 
of the Banach algebras are independent of $p$.}


\section{Introduction}

Let $G$ be a locally compact group. Denote $\underline E G$ its classifying space for proper actions, and $K^G_i(\underline E G)$ the equivariant K-homology groups. The Baum-Connes conjecture claims that the associated Baum-Connes assembly map $\mu: K^G_i(\underline E G)\to K_i(C^*_r(G))$ is an isomorphism between abelian groups, where $K_i(C^*_r(G))$ is the $i$-th K-theory group of the reduced $C^*$ algebra $C^*_r(G)$ of $G.$ The conjecture has been verified for a large class of groups, including groups with Haagerup property \cite{higson-kasparov}, hyperbolic groups \cite{mineyev-yu,lafforgue02}, reductive Lie groups over local fields and
cocompact lattices in $SL_3$ over a local field \cite{lafforgue02}.

If one replaces $C^*_r(G)$ by the Banach algebra $L^1(G)$ of integrable functions on $G$, this is the so-called Bost's conjecture, which has been proved for a much larger class of groups, including all lattices in a reductive Lie group over a local field \cite{lafforgue02}.

Motivated by a theorem due to the second named author that hyperbolic groups act isometrically properly on some $L^p$ space \cite{yu05,nica13}, G.~Kasparov and he introduced an $L^p$ version of Baum-Connes assembly map, where one replaces $C^*_r(G)$ by the Banach algebra of convolution operators on the space of $p$-integrable functions $L^p(G)$ on $G$ \cite{kasparov-yu-draft}. This algebra is denoted by $B^p_r(G)$ in this article, which in more precise terms is the Banach algebra obtained by completing the convolution algebra of compactly supported continuous functions $C_c(G)$ with respect to the operator norm on $L^p(G)$.

These algebras of operators in $L^p$-spaces also appear in the work of N.~C.~Phillips \cite{phillips13} on the study on $L^p$ cross-product and Cuntz algebras.

The aim of this work is to investigate the relations among the K-theory groups of these different group Banach algebras $B^p_r(G)$.

\begin{thm}\label{duality}
For a locally compact group $G,$
the K-theory groups of $B^p_r(G)$ and $B^q_r(G)$
are canonically isomorphic, where $p\in [1,\infty]$ and 
$q$ is the dual number of $p,$ namely
$1/q+1/p=1.$
\end{thm}


To compare K-theory groups of algebras $B^p_r(G)$ for different $p,$ a natural way is to consider a subalgebra $B^{p,*}_r(G)$
in $B^p_r(G)$ which is closed under involution. This is the involutive Banach algebra obtained by completing $C_c(G)$ with respect to the norm 
$$\|f\|_{B^{p,*}_r(G)}=\max\{\|f\|_{B(L^p(G))},\|f^*\|_{B(L^p(G))}\}$$
where $f^*(g)=\overline{ f(g^{-1})}$.
Under the condition that the group has the property of rapid decay (property RD) \cite{jolissaint90,chatterji-pittet-saloff-coste,valette02,sapir}, we are able to show that this procedure does not lose K-theoretic information. 

\begin{thm}\label{subalgebra}
Let $p\in [1,\infty].$ 
Let $G$ be a locally compact group with property RD. 
The K-theory groups of $B^p_r(G)$ and $B^{p,*}_r(G)$ are canonically isomorphic.
\end{thm}

In fact, the statement holds 
under a weaker condition that the group has a Banach version of property RD, which we call property $(RD)_q$ where $1/q+1/p=1$ (Theorem~\ref{rd-rdq} and Theorem~\ref{*-independent-p}).
We prove this by showing that there exists a subalgebra which is closed under holomorphic functional calculus in both $B^p_r(G)$ and $B^{p,*}_r(G)$ (Proposition~\ref{St}).

In the case when the group acts properly on an $L^p$ space, 
this phenomena of $*$-independence is also verified in \cite{kasparov-yu-draft}.
We conjecture that it is a general phenomenon for groups.

The advantage of $B^{p,*}_r(G)$ is that it allows to use
interpolation to give a canonical morphism 
$$i_{p',p}:B^{p',*}_r(G)\to B^{p,*}_r(G)$$ 
for $2\leq p<p'\leq \infty.$
Letting $p$ vary from $1$ to $2,$ we obtain a family of continuous inclusions 
$$L^1(G)\xrightarrow{i_{1,p}} B^{p,*}_r(G)\xrightarrow{i_{p,2}} C^*_r(G).$$




When the group is discrete,
we show the following semi-continuity result for $p\to 2$.

\begin{thm}\label{semi-cont}
Let $G$ be a finitely generated group.
$K_*(C^*_r(G))$ is the inductive limit of the system $\{K_*(B^{p,*}_r(G)),p>2\}$.
Namely,
$$K_*(C^*_r(G))=\bigcup_{p>2}i_{p,2}^*(K_*(B^{p,*}_r(G))).$$
\end{thm}

The following corollary gives us an alternative criterion for the K-theory of $C^*_r(G)$: it is either equal to that of $B^{p,*}_r(G)$ for some $p>2,$ or it is dramatically larger than that of every $B^{p,*}_r(G),p>2.$
\begin{cor}
Let $G$ be a finitely generated group. Then 
there exists $p>2$ such that 
$$i^*_{p,2}(K_*(B^{p,*}_r(G)))= K_*(C^*_r(G)),$$
unless for every $p>2$ (in particular $K_*(l^1(G))$) $i^*_{p,2}(K_*(B^{p,*}_r(G)))\subset K_*(C^*_r(G))$ is of infinite index.
\end{cor}


For some large class of groups listed below, via a Banach version of property RD that we call property $(RD)_q$ (see Section 3),
and unconditional completion \cite{lafforgue02}, we are able to show that their K-theory groups indeed do not depend on the parameter $p\in [1,\infty].$

\begin{thm}\label{p-independent}
Let $G$ be either 
\\- a semi-simple Lie group over a local field, or a cocompact lattice in it with property RD, or
\\- a hyperbolic group, or
\\- a finitely generated group of polynomial growth. 
\\
$K_*(B^p_r(G))$ are isomorphic for $p\in [1,\infty]$.
\end{thm}


{\bf Acknowledgement: } This work is supported in part by NSF and CNSF 11420101001. The authors would like to thank Vincent Lafforgue for telling us about Theorem~\ref{rd-rdq} and other helpful comments. They would also like to thank Gilles Pisier for sharing his insight in $L^p$ convolution operators, and providing an interpolation proof of Theorem~\ref{rd-rdq}.

\section{The Banach algebras $B^{p}_r(G)$ and $B^{p,*}_r(G)$}

In this section, we introduce the $L^p$ convolution algebra $B^p_r(G)$ and its involutive counterpart $B^{p,*}_r(G)$, and investigate their basic properties.

Let $G$ be a locally compact group, $p\in [1,\infty].$ 
Denote $B^p_r(G)$ the Banach algebra obtained by
completing $C_c(G)$ 
with the operator norm $\|f\|_{B(L^p(G))},f\in C_c(G)$.
Denote $B^{p,*}_r(G)$ the involutive version of $B^p_r(G)$, namely the completion of $C_c(G)$ with respect to the norm 
$$\|f\|_{B^{p,*}_r(G)}=\max\{\|f\|_{B(L^p(G))},\|f^*\|_{B(L^p(G))}\},$$
where $f^{*}(x)=\overline{f(x^{-1})}.$

Except for the abelian case, these two algebras are usually not identical.

\begin{lem}\label{sg}
If $G$ is a countable discrete group and $H\subset G$ a subgroup, then for any $f\in \C(G)$ supported in $H$ we have $\|f\|_{B^p_r(G)}=\|f\|_{B^p_r(H)}$.
\end{lem}
\qed

\begin{prop}
If $G$ is a non-elementary hyperbolic group, or a non-abelian linear group, then 
$$B^{p,*}_r(G)\neq B^p_r(G)$$
for any $p\in (1,2)\cup (2,\infty)$.
\end{prop}
{\bf Proof}: 
By Tit's alternative, a non-abelian linear group contains a non-abelian free subgroup. The same holds as well for a non-elementary hyperbolic group.
The statement then follows from Lemma~\ref{sg} and the main Theorem in \cite{pytlik}.

The argument obviously works for any discrete group containing a non-abelian free group as a subgroup.
\qed


\begin{prop}
There exists an amenable group $G$ such that
$$B^{p,*}_r(G)\neq B^p_r(G)$$
for $p=4.$
\end{prop}

{\bf Proof}: 
The construction of $G$ is based on results in \cite{oberlin}. We thank G.~Pisier for pointing us to \cite{oberlin} and generously sharing his observations.

\begin{lemma}\label{oberlin}\cite{oberlin}
There exists a finite group $G_0$ and a function $f\in~\C (G_0)$ such that
$$\|f^*\|_{B(l^4(G_0))}>\|f\|_{B(l^4(G_0))}.$$
\end{lemma}
\qed


In what follows, we construct the amenable group $G$ as in the statement of the proposition,
and functions $f_n\in\C (G)$ such that
$$\|f_n\|_{B(l^4(G))}=1$$
and
$$\|f_n^*\|_{B(l^4(G))}\to \infty.$$

Define 
$$G_n:=\prod_{i=1}^nG_0,$$
where $G_0$ is as in Lemma~\ref{oberlin},
and
$$G=\{(g_i\in G_0),\tx{ only finitely many }g_i\tx{ are not the neutral element}\}\subset \prod_{i=1}^\infty G_0.$$
The following map
$$G_n\to G$$
$$(g_1,...,g_n)\mapsto (g_1,...,g_n,e,...,e,...)$$
is an embedding of groups. We have $G_n\subset G_{n+1}$.
Define $f_n':G_n\to \C$ by
$$f_n'(g_1,...,g_n)=f(g_1)...f(g_n)/M^n$$
where 
$$M=\|f\|_{B(l^4(G_0))}.$$
Extend the definition from $G_n$ to $G$ by zero outside of $G_n$, we get a finitely supported function $f_n\in \C (G).$

\begin{lemma}
Let $H_1,H_2,$ be discrete groups, and $\phi_i\in \C (H_i).$ The tensor product $\phi=\phi_1\otimes\phi_2$
(namely $\phi(h_1,h_2)=\phi(h_1)\phi(h_2)$)
is finitely supported on $H=H_1\times H_2,$
and we have for all $p\in [1,\infty]$
$$\|\phi\|_{B(l^p(H))}=\|\phi_1\|_{B(l^p(H_1))}\|\phi_2\|_{B(l^p(H_2))}.$$
\end{lemma}
\qed

Let $G_n'$ be the subgroup in $G$ such that $G=G_n\times G_n'$. Since $f_n=f_n'\otimes \delta_{e_{G_n'}}$, we have
$\|f_n\|_{B(l^4(G))}= 1$
and 
$$\|f_n^*\|_{B(l^4(G))}=(\|f^*\|_{B(l^4(G_0))}/M)^n\to \infty.$$

This terminates the proof.
\qed

\begin{prop}
Let $G$ be a locally compact group.
For 
$p<p'<~\infty$,
the identity map on $C_c(G)$ extends to a continuous (contractive) injective morphism of Banach algebras
$$i_{p',p}:B^{p',*}_r(G)\to B^{p,*}_r(G).$$
\end{prop}

{\bf Remark. } $i_{p',p}$ fails to be surjective in general. Indeed, it is shown in \cite{pytlik} that there exists a function $f$ on a non-abelian free group $G$ such that $\|f\|_{C^*_r(G)}$ is finite and $\|f\|_{B(l^{p'}(G))}$ is infinite for any $p'\neq 2,$ which means that $i_{p',2}$ is never surjective on $G.$
\\

{\bf Proof. }
Denote $q,q'$ the duals of $p,p'$ respectively
$:1/p+1/q=1,$ $1/p'+1/q'=1.$ Let $\theta\in [0,1]$ such that $1/p=\theta /p'+(1-\theta)/q',$ as a consequence $1/q=(1-\theta)/p'+\theta/q'.$
By complex interpolation we have
$$\|f\|_{B(L^p(G))}\leq \|f\|_{B(L^{p'}(G))}^\theta\|f\|_{B(L^{q'}(G))}^{1-\theta},$$
and 
$$\|f\|_{B(L^q(G))}\leq \|f\|_{B(L^{p'}(G))}^{1-\theta}\|f\|_{B(L^{q'}(G))}^{\theta}.$$
Consequently, we have
$$\|f\|_{B^{p,*}_r(G)}\leq \|f\|_{B^{p',*}_r(G)}.$$
So $i_{p',p}$ is continuous and contractive. 

Now prove that $i_{p',p}$ is injective. Suppose $F\in {B^{p',*}_r(G)}$ such that $i_{p',p}(F)=0.$ Prove that $F(\xi)=0$ for any $\xi\in L^{p'}(G)$ and $F(\eta)=~0,$ $\forall \eta\in L^{q'}(G)$. 
Since $C_c(G)$ is a dense subset in both $L^{p'}(G)$ and $L^{q'}(G)$,
by continuity of $F,$ it suffices to prove that 
$\pair{\phi,F(f)}=0,$ $\forall f,\phi\in C_c(G).$
 Let $F_n\in C_c(G)$ such that $F_n\to F$ in $B^{p',*}_r(G)$. By continuity of $i_{p',p}$ we have $F_n\to 0$ in $B^{p,*}_r(G)$. Thus $\pair{\phi,F(f)}=\lim_n \pair{\phi,F_n(f)}=0.$
\qed

\section{Semi-continuity of $K_*(B^{p,*}_r)$}

In this section, we investigate K-theoretic properties of $B^{p}_r(G)$ and $B^{p,*}_r(G)$. 
We first show Theorem~\ref{duality}.

{\bf Proof of Theorem~\ref{duality}:}
Observe that 
$$*:B^p_r(G)\to B^q_r(G),f\mapsto [g\mapsto \overline{f(g^{-1})}]$$
is an isometric anti-isomorphism between Banach algebras. Thus it sends idempotents onto idempotents, invertibles onto invertibles, and preserves equivalence relation.
\qed
\\

Next we prove Theorem~\ref{semi-cont}, which is merely a special case of the following theorem.

\begin{thm}
Let $G$ be a finitely generated group.
For $j=0,1$,
$i_{p',p}^*(K_j(B^{p',*}_r(G)))$ is a subgroup in $K_j(B^{p,*}_r(G))$, and we have
$$K_j(B^{p,*}_r(G))=\bigcup_{p'>p}i_{p',p}^*(K_j(B^{p',*}_r(G))).~~~~(*)$$
\end{thm}


{\bf Proof. }
The following argument is inspired by \cite{lafforgue02}.


Let $\alpha\in [0,1]$ such that $1/q'=(1-\alpha)+\alpha/q$ (consequently $1/p'=\alpha/p$). We show that for $F\in M_n(\C(G))$
$$\|F\|_{M_n(B^{p',*}_r(G))}\leq e^{\frac{\lambda(1-\alpha)}{q}m_F }\|F\|_{M_n(B^{p,*}_r(G))},~~~~(**)$$
where $\lambda >0 $ is such that $|B_m(0)|\leq e^{\lambda m}$ and $m_F$ is the smallest $m$ such that all $F_{i,j}\in \C(G),1\leq i,j\leq n,$ are supported in $B_m(0)\subset G$.
It suffices to show inequality $(**)$ for $n=1.$
Let $f\in \C(G)$ with support in $B_m(0).$ By complex interpolation, we have
$$\|f\|_{B(\ell^{q'}(G))}\leq \|f\|_{\ell^1(G)}^{1-\alpha} \|f\|_{B(\ell^{q}(G))}^{\alpha},$$
$$\|f\|_{B(\ell^{p'}(G))}\leq \|f\|_{B(\ell^\infty(G))}^{1-\alpha} \|f\|_{B(\ell^{p}(G))}^{\alpha}\leq
\|f\|_{\ell^1(G)}^{1-\alpha} \|f\|_{B(\ell^{p}(G))}^{\alpha}.$$
We have for any $s>1,$ $\|f\|_{\ell^s(G)}\leq \|f\|_{B(\ell^s(G))},$ and by Holder's inequality, $\|f\|_{\ell^1(G)}\leq e^{\lambda m/t}\|f\|_{\ell^s(G)}$ where $1/t+1/s=1.$ 
Apply $s=p,q$ and combine with the previous two inequalities, we obtain the desired $(**)$ (for $n=1$).

Prove $(*)$ for $j=0.$
It suffices to show that, for $n\in\N$ and $e\in Idem(M_n(B^{p,*}_r(G)))$, there exist $p'>p$ and $e'\in Idem(M_n(B^{p',*}_r(G)))$ such that 
$$\|e-e'\|_{M_n(B^{p,*}_r(G))}<1/\|1-2e\|_{M_n(B^{p,*}_r(G))},$$
where $\|F\|_{M_n(B)}$ is defined to be $\max_{i,j}\|F_{i,j}\|_{B}$ for a Banach algebra $B$.
Let $e\in Idem (M_n(B^{p,*}_r(G)))$. By density, for any $\eps>0$ there exists $f=f_\eps \in M_n(\C(G))$ such that $\|e-f\|_{M_n(B^{p,*}_r(G))}<\eps.$ Denote $M=\|e\|_{M_n(B^{p,*}_r(G))}$.
Let $m\in\N$ such that $B_m(0)$ contains the support of all $(f^2-f)_{i,j}$. 
Since $f^2-f=(f-e)f+(e-1)(f-e)$ we have
$\|f^2-f\|_{M_n(B^{p,*}_r(G))}< \eps(2M+\eps+1).$
Consider those $\eps>0$ such that $\eps(2M+\eps+1)<1/4.$
 By inequality $(**)$ there exists $p'>p$ (depending on $m$ and $M$) such that ($1-\alpha>0$ is small enough and consequently) 
$$\|f^2-f\|_{M_n(B^{p',*}_r(G))}< 1/4.$$
Let $D_0,D_1$ be the open disk of radius $1/2$ and centered at $0,1$ respectively on the complex plane, and
denote by $\varphi$ the analytic function on $U=D_0\sqcup D_1$ that sends $D_0$ to $0$ and $D_1$ to $1.$ We have that $Spec_{M_n(B^{p',*}_r(G))}(f)\subset U,$ and $\varphi (f)$ is an idempotent in $M_n(B^{p',*}_r(G)).$ By continuity of the inverse map on $B^{p,*}_r(G),$ we have that $\varphi (f_\eps)\to \varphi(e)=e$ when $\eps \to 0.$
Take $\eps>0$ small enough such that $\|\varphi(f_\eps)-e\|_{M_n(B^{p,*}_r(G))}<1/\|1-e\|_{M_n(B^{p,*}_r(G))}$ and we are done.

Now Prove $(*)$ for $j=1.$ Let $F\in GL_n(B^{p,*}_r(G))$ be an invertible element. Choose $f\in M_n(\C(G))$ close to $F$ and $\phi\in M_n(\C(G))$ close to $F^{-1}$ such that $\|\phi f-1\|_{M_n(B^{p,*}_r(G))}\leq 1/3$ and $f$ lies in the same connected component as $F$ in $GL_n(B^{p,*}_r(G))$.
By inequality $(**)$ there exists $p'>p$ such that $\|\phi f-1\|_{M_n(B^{p',*}_r(G))}\leq 1/2.$ Thus $\phi f$ is invertible in $M_n(B_r^{*,p'}(G))$ and consequently $f^{-1}=(\phi f)^{-1}\phi\in M_n(B_r^{*,p'}(G)).$ 


\qed
\\

\section{A Banach version of property $RD$}

In this section, we define a Banach version of property RD, and prove Theorem~\ref{p-independent} in the introduction.

\begin{defi} Let $q\in [1,\infty]$ and $G$ be a locally compact group.
Say that $G$ has $(RD)_q$ 
(with respect to a measurable length function $l$),
if there exists a polynomial $P$ such that for any continuous function $f$ with support in $B_n(e)$, we have
$$\|f\|_{B(L^q(G))}\leq P(n)\|f\|_{L^q(G)}.$$
\end{defi}

It is obvious that every group has $(RD)_1.$ 

\begin{prop}
Let $G$ be a locally compact group. If $G$ is of polynomial growth with respect to some length function $l,$ then
it has $(RD)_q$ for all $q\in [1,\infty].$
\end{prop}

{\bf Proof. }
By Holder's inequality, we have for any continuous function $f$ supported in $B_n(e)$
$$\|f\|_{B(L^q(G))}\leq \|f\|_{L^1(G)}\leq |B_n(e)|^{1/p}\|f\|_{L^q(G)}.$$
A polynomial bound on $|B_n(e)|$ clearly yields $(RD)_q$ for $G.$
\qed
\\

The following proposition suggests that the notion of $(RD)_q$ is 
more interesting
when $q\in (1,2].$

\begin{prop} \label{p>2}
A countable discrete group having $(RD)_p$ for some $p\in (2,\infty]$ with respect to a length function $l$
is of polynomial growth in $l$.
\end{prop}


{\bf Proof. }
Let $q\in [1,2)$ be the dual number of $p.$
Let $P$ be the polynomial for $(RD)_p$ of $G$. We have
for $f\in\C (G)$ supported in $B_n(e)$
$$\|f\|_{l^q(G)}\leq \|f\|_{B(l^q(G))}=\|f^*\|_{B(l^p(G))}$$
$$\leq P(n)\|f^*\|_{l^p(G)}=P(n)\|f\|_{l^p(G)}.$$
Take $f_n=\chi_{B_n}$, we have 
$$|B_n|^{1/q}\leq P(n)|B_n|^{1/p}.$$
Since
$1/q-1/p>0,$ it is immediate that the group has polynomial growth.
\qed
\\


The following result is due to V.~Lafforgue. His original proof involves a combinatorial characterization of RD. Here we give a proof in the discrete case based on an idea using Mazur map, and another proof due to G.~Pisier based on complex interpolation. We think that these two proofs are different and have their own merits, and include them both in the article. B.~Nica informed us that he also has a proof of this result.

\begin{thm} (V.~Lafforgue)\label{rd-rdq}
If $G$ is a locally compact group with property $(RD)_q$ for some $q>1$, 
then it has $(RD)_{q'}$ for any $q'\in (1,q).$
In particular, RD implies $(RD)_q$ for $q\in (1,2).$
\end{thm}

{\bf First proof}
(in the discrete case). 
Suppose that $G$ is a discrete group having $(RD)_q$ with respect to a polynomial $P.$
Let $\phi\in \C(G)$.
Set 
$$\phi_\alpha (g):=|\phi(g)|^\alpha, \alpha=q'/q<1.$$
We then have
$\|\phi\|_q^q=\|\phi_\alpha\|^{q'}_{q'}.$
Since
$a_1^\alpha+...+a_n^\alpha\geq (a_1+...+a_n)^\alpha$ for $a_i\geq 0$ and $\alpha <1$,
we have
$$f_\alpha\phi_\alpha(g)\geq (|f||\phi|)_\alpha(g).$$
Therefore,
$$\||f||\phi|\|^{q'}_{q'}=\|(|f||\phi|)_\alpha\|_q^q\leq \|f_\alpha\phi_\alpha\|_q^q$$
$$\leq P(n)^q\|f_\alpha\|_q^q\|\phi_\alpha\|_q^q=P(n)^q\|f\|_{q'}^{q'}\|\phi\|_{q'}^{q'}.$$
This completes the proof since the left hand side is $\geq \|f\phi\|_{q'}^{q'}.$ 
\qed
\\

{\bf Second proof} (G.~Pisier). 
Let $\theta\in (0,1)$ such that
$1/q'=(1-\theta)+\theta/q$. Let $f\in C_c(G)$ be a function supported in the ball of radius $n$ and define
$$F_z(g)=|f(g)|^{(1-z+z/q)q'}\phi(g)$$
where $\phi$ is the phase function of $f,$ namely $\phi(g)=f(g)/|f(g)|$ whenever $f(g)\neq 0$ and zero otherwise. We have obvious relations
$$F_\theta(g)=f(g),
F_0(g)=|f(g)|^{q'}\phi(g),
F_1(g)=|f(g)|^{q'/q}\phi(g),$$ 
and
$\|F_0\|_1=\|F_1\|_{q}^q=\|F_\theta\|_{q'}^{q'}.$

Let $x,y\in C_c(G)$ be arbitrary functions. Define $X_z,Y_z$ for $x,y$ by the same formula as $F_z$ for $f.$
Now the function $\{z\in \C,0\leq Re(z)\leq 1\}\to \C,z\mapsto \pair{X_z,F_z Y_z}$ is analytic on the interior and continuous on the boundary.
By Hadamard three-lines lemma, we have
$$ |\pair{x,fy}|=
|\pair{X_\theta,F_\theta Y_\theta}|\leq \sup_{t\in \R}\pair{X_{it},F_{it} Y_{it}}^{1-\theta}
\sup_{t\in \R}\pair{X_{1+it},F_{1+it} Y_{1+it}}^\theta$$
$$\leq (\|x\|_{q'}\|f\|_{q'}\|y\|_{q'})^{(1-\theta)q'}P(n)^\theta(\|x\|_{q'}\|f\|_{q'}\|y\|_{q'})^{\theta q'/q}$$
$$=P(n)^\theta\|x\|_{q'}\|f\|_{q'}\|y\|_{q'},$$
where $P$ is the polynomial function in the definition $(RD)_q$.
This implies $\|f\|_{B(L^{q'}(G))}\leq P(n)^{\theta}\|f\|_{L^{q'}(G)}$.
\qed
\\

Analogous to the $L^2$ case, we have the following result for amenable groups with rapid decay.

\begin{prop}
Let $G$ be a compactly generated amenable group. If $G$ has $(RD)_q$ for some $q\in (1,2],$ then $G$ is of polynomial growth.
\end{prop} 

{\bf Proof. }
The following argument is well-known.

\begin{lem*}
Let $G$ be amenable and $f\in C_c (G)$ a non negative function. Then $$\|f\|_{B(L^p(G))}=\|f\|_{L^1(G)}$$ for any $p\in [1,\infty]$.
\end{lem*}

{\bf Proof. }
The statement follows from applying Folner sets.
\qed

Now take $f=\chi_{B_n}$. By similar arguments as in the proof of Proposition~\ref{p>2} one sees that $G$ is of polynomial growth.
\qed
\\



\begin{prop}\label{St}
Let $p\in [1,\infty]$ and $q$ its dual number.
Let $G$ be a locally compact group with property $(RD)_q$ with respect to a continuous length function $L.$
Then 
for sufficiently large $t>0,$
the space $S^t_q(G)$ of elements $f\in L^q(G)$ such that 
$$\|f\|_{S^t_q}:=\|g\mapsto (1+L(g))^tf(g)\|_{L^q(G)}<\infty$$
is a Banach algebra for the norm $\|\cdot\|_{S^t_q}$.
It is
contained in $B^{p,*}_r(G)\subset B^p_r(G),B^q_r(G)$, 
and stable under holomorphic functional calculus in each of these three algebras.
\end{prop}

{\bf Proof. }
Here we generalize the argument in \cite{lafforgue00} to the case of locally compact groups and $L^p$ spaces.

Suppose $G$ has $(RD)_q$ with polynomial $n\mapsto Cn^D$ for some $C,D>0.$
First prove containment, namely
$$\|f\|_{B(L^q(G))}\leq K\|f\|_{S^t_q},t\geq D+1$$
(by duality the inequality for $\|f\|_{B(L^p(G))}$ follows from this one).
For this let
$S_n:=B_n(e)\backslash B_{n-1}(e),f_n:=f1_{S_n}$
$$\|f\|_{BL^q}\leq \sum_n \|f_n\|_{BL^q}\leq \sum_n Cn^D\|f_n\|_{L^q}$$
by Holder's inequalty
$$\leq C(\sum n^{-p})^{1/p}(\sum n^{(D+1)q}\|f_n\|^q_{L^q})^{1/q}
=K\|f\|_{S^{D+1}_q}.$$

Now we prove that it is an algebra.
By the inequality $(x+y)^s\leq 2^{s-1}(x^s+y^s),s\geq 1$ and the triangular inequality for $L$, we have
for non-negative $f_1,f_2\in C_c(G),x\in G$
$$f_1*f_2(x)(1+L(x))^s\leq 2^{s-1} \Big(f_1*\big (f_2(1+L)^s\big)(x)+\big (f_1(1+L)^s\big)*f_2(x)\Big).$$
Therefore,
$$\|f_1*f_2\|_{S^s_q}^q\leq 2^{q(s-1)+q-1}(\|f_1*\big (f_2(1+L)^s\big)\|_{L^q}^q+\|\big (f_1(1+L)^s\big)*f_2\|_{L^q}^q)$$
$$\leq 2^{qs-1}K^q(\|f_1\|_{S^{D+1}_q}^q\|f_2\|^q_{S^s_q}+\|f_2\|_{S^{D+1}_q}^q\|f_1\|^q_{S^s_q})\leq 2^{qs}K^q\|f_1\|_{S^s_q}^q\|f_2\|_{S^s_q}^q,$$
namely
$$\|f_1*f_2\|_{S^s_q}\leq 2^sK\|f_1\|_{S^s_q}\|f_2\|_{S^s_q}$$
for $s\geq D+1.$
The inequality still holds for $f_1,f_2\in C_c(G)$ without the non-negativity assumption by triangular inequality.

Now prove that it is stable under holomorphic functional calculus in both $B^q_r(G)$ and $B^p_r(G)$ (and consequently it is so in $B^{p,*}_r(G)$ as well). In fact, for $B^q_r(G)$ it suffices to show 
for $f\in C_c(G)$
$$\lim_n\|f^n\|^{1/n}_{S^t_q}=\lim_n\|f^n\|^{1/n}_{B(L^q(G))}.$$
First of all, by an argument as before we have
$$\|f^{n+1}\|_{S^s_q}^q\leq 2^{qs-1}(\|f\|_{B(L^q(G))}^q\|f^n\|^q_{S^s_q}+\|f^n\|_{B(L^q(G))}^q\|f\|^q_{S^s_q}).$$
By an inductive argument we have $\forall n$
$$\|f^n\|_{S^s_q}\leq 2^{ns}\|f\|_{S^s_q}\|f\|^{n-1}_{B(L^q(G))},$$
which implies 
$$\lim_n\|f^n\|_{S^s_q}^{1/n}\leq 2^s\|f\|_{B(L^q(G))}.$$
Replacing $f$ by $f^p$ and let $p\to \infty$ we have proved $\leq$ part in the inequality. The other part $\geq$ follows from the inequality concerning the containment $S^s_q(G)\subset B^{q}_r(G)$ as in the beginning of this proof. The proof for $B^p_r(G)$ is similar and so is omitted.
\qed
\\

Let $G$ be as before a locally compact group and $L:G\to \R$ a continuous length function. Denote $C_u(G)$ the translation algebra, namely $C_c(G,L^\infty(G))$, and $B^{p}_u(G),p\in [1,\infty]$ its operator completion on $L^p(G)$, and $B^{p,*}_u(G)$ its involutive counterpart.

Let $\delta$ be the densely defined derivation on $B^{p}_u(G)$ defined by the formula 
$$\delta (T)=[l,T]=lT-Tl,$$
where by abuse of notation, $l$ is the densely defined operator on $L^p(G)$ by point-wise multiplication $l:f\mapsto [x\mapsto l(x)f(x)]$. $\delta$ is also a densely defined derivation on $B^{p,*}_u(G)$ by the same formula.

The following proposition is well-known when $p=2$, yet we cannot find a proof in the literature.

\begin{prop}
$\delta$ is a closed derivation on both $B^{p}_u(G)$ and $B^{p,*}_u(G)$.
\end{prop}

{\bf Proof. }
Following 
\cite{bost}
we say that that a one-parameter group of automorphisms $a_*:\R\to Aut(B)$ on a Banach algebra is strongly continuous if $[t\mapsto a_tb]\in C(\R;B)$ is a continuous mapping for any $b\in B.$ Its generator is denoted by $\delta_a.$


\begin{lem*} (\cite{bost} Lemma~4.2.1) 
Let $a_*:\R\to Aut(B)$  be a strongly continuous as above. Then $x\in Dom(\delta_a),y=\delta_ax\in B$ if and only if 
$$a_tx=x+\int_0^ta_syds,\forall t\in \R.$$
As a consequence, $\delta_a$ is a closed operator.
\end{lem*}
\qed


Recall that $l:G\to \R$ is 
a continuous length function.
Define the following algebraic automorphism
$$a_t:B^p_u(G)\to B^p_u(G),T\mapsto e^{ilt}Te^{-ilt},t\in \R$$
where $e^{itl}$ acts on $L^p(G)$ unitarily by point-wise multiplication of the function $x\mapsto e^{itl(x)}$.
It is clear that $\delta_a=i\delta.$



By a density argument, it suffices to prove continuity on $C_u(G).$
First
$$|e^{it}-1|^2=(cost-1)^2+(sint)^2\leq C^2t^2$$
for some $C>0.$
Therefore,
$$(a_te_g-e_g)v(x)=(e^{itl(x)-itl(g^{-1}x)}-1)v(g^{-1}x),$$
and
$$\|(a_te_g-e_g)v\|_p\leq Ctl(g)\|v\|_p$$
This implies that for $T=\int e_gf_gdg$ with finite propagation
$$\|a_t(T)-T\|_{B(L^p(G))}\leq Ct\int l(g)\|f_g\|_\infty dg\to 0,t\to 0.$$
Therefore, $\delta$ is a closed derivation on both $B_u^{p}(G)$ and $B_u^{p,*}(G)$.
\qed
\\


\begin{prop}\label{S-infinity}
Let $p\in [1,\infty]$ and $q$ be its dual number. When $G$ is discrete,
$S^\infty_q(G):=\cap_{t>0}S^t_q(G)$ is the space of smooth vectors in $B^p_r(G),B^q_r(G)$ and $B^{p,*}_r(G)$ with respect to $\delta$.
\end{prop}


{\bf Proof. }
The following argument is well-known. We include it here for the reader's convenience.

First of all, we have
$$(\delta^k(T)\xi) (g)=\sum_h T_{g,h}\xi(h)(l(g)-l(h))^k,$$
which implies that for $g\in G$
$$|(\delta^k(f)\xi)(g)|\leq ((l^k|f|)*|\xi|)(g)$$
and therefore, by $(RD)_q$
$$\|\delta^k f\|_{B(l^r(G))}\leq \|(l^k|f|)\|_{B(l^r(G))}\leq \|f(1+l)^{k+s}\|_{l^q(G)},r=p,q$$
for sufficiently large $s>0$.
Thus 
$S^\infty_q$
is contained in the intersection of domains of $\delta^k$ in $B^p_r(G)$, $B^q_r(G)$ and also $B^{p,*}_r(G)$.

For the other inclusion, notice that 
$$\|f l^k\|_{l^q(G)}=\|\delta^k f(\delta_e)\|_{l^q(G)}\leq \|\delta^kf\|_{B(l^q(G))},$$
$$\|f l^k\|_{l^q(G)}=\|f^* l^k\|_{l^q(G)}=
\|\delta^k f^*(\delta_e)\|_{l^q(G)}\leq \|\delta^kf\|_{B(l^p(G))},$$
where $\delta_e$ denotes the Dirac function at $e\in G.$
 
In summary, $S^\infty_q(G)$ is exactly the space of smooth vectors in $B^p_r(G)$, $B^q_r(G)$ and $B^{p,*}_r(G)$.
\qed
\\

Theorem~\ref{subalgebra} in the introduction is a special case of the following statement.

\begin{cor}\label{*-independent-p}
Let $q_o\in [1,\infty],$ $q\in [1,q_o]$, $p$ be the dual of $q.$
Let $G$ be a locally compact group with property $(RD)_{q_o}.$  
Then the canonical algebraic morphism $B^{p,*}_r(G)\to B^{p}_r(G)$ induces isomorphism between their K-theory groups. So does $B^{p,*}_r(G)\to B^{q}_r(G)$.
\end{cor}

{\bf Proof. }
The statement follows from Theorem~\ref{rd-rdq} and Proposition~\ref{St}.
\qed
\\

Theorem~\ref{p-independent} is a special case of the following Corollary.

\begin{cor}
Let $G$ be a locally compact group in Lafforgue's class C' \cite{lafforgue02} having property $(RD)_q,q\in [1,2].$ 
$i_{p',p}:B^{p',*}_r(G)\to B^{p,*}_r(G)$
induces isomorphism in K-theory for any $2\leq p<p'\leq+\infty,$ where $p$ is the dual number of $q.$
\end{cor}

{\bf Proof. }
Apply the assembly maps for unconditional completions \cite{lafforgue02} to $S^t_q(G)$ and $S^t_{q'}(G)$ for sufficiently large $t>0$,
where $q'$ is dual to $p'.$
The claim follows from Proposition~\ref{St}, 
the injectivity \cite{kasparov88,kasparov-skandalis91,kasparov-skandalis03} and surjectivity \cite{mineyev-yu,lafforgue02} of the assembly maps \cite{lafforgue02}.
\qed

\section{Open problems}

In this last section, we list several 
interesting 
open problems.

\begin{enumerate}
\item Is the algebra $B^{p}_r(G)$ non-involutive for every (discrete) non-amenable group? Namely, is it true that on a non-amenable group, there exists a function $f:G\to \C$ that acts as a bounded operator by convolution on $L^p(G)$, but its involution $f^*$ does not?
\item Does the canonical morphism $B^{p,*}_r(G)\to B^p_r(G),p\geq 2$, always induce an isomorphism in K-theory? We already know that it is true for groups with $(RD)_q$ where $1/p+1/q=1$, and for groups acting properly isometrically on $L^p$ spaces.
\item Is it true that $K_*(B_r^p(G))$ and $K_*(B_r^{p,*}(G))$ are independent of $p$ for any locally compact group $G$?
\item Is there a locally compact group with property $RD$ but fails $(RD)_q$ for some $q\in (1,2)$?
\item Does semi-continuity (Theorem~\ref{semi-cont}) hold for a locally compact group? What about crossed product? See \cite{chung} for a discussion on $L^p$ cross product.
\end{enumerate}

\bibliographystyle{plain}
\bibliography{/Users/ben/Dropbox/my-tex/bibtex_lib}

\begin{thebibliography}{10}

\bibitem{bost}
J.-B. Bost.
\newblock {Principe d'Oka, K-th{\'e}orie et syst{\`e}mes dynamiques non
  commutatifs}.
\newblock {\em Inventiones mathematicae}, 101(2):261--334, 1990.

\bibitem{chatterji-pittet-saloff-coste}
I.~Chatterji, C.~Pittet, and L.~Saloff-Coste.
\newblock {Connected Lie groups and property RD}.
\newblock {\em Duke Math. J.}, 137(3):511--536, 04 2007.

\bibitem{chung}
Y.-C. Chung.
\newblock {Dynamic asymptotic dimension and K-theory of Banach crossed product
  algebras}.
\newblock {\em https://arxiv.org/abs/1611.09000}.

\bibitem{kasparov-yu-draft}
G.~Yu. G.~Kasparov.
\newblock {The Baum-Connes conjecture and group The Baum-Connes conjecture and
  group actions on lp-spaces}.
\newblock {\em In preparation}.

\bibitem{higson-kasparov}
Nigel Higson and Gennadi Kasparov.
\newblock {E-theory and KK-theory for groups which act properly and
  isometrically on Hilbert space}.
\newblock {\em Invent. math.}, 144:23--74, 04 2001.

\bibitem{jolissaint90}
Paul Jolissaint.
\newblock {Rapidly Decreasing Functions in Reduced C*-Algebras of Groups}.
\newblock {\em Transactions of the American Mathematical Society},
  317(1):167--196, 1990.

\bibitem{kasparov-skandalis91}
G.~G. Kasparov and G.~Skandalis.
\newblock {Groups acting on buildings, operator K-theory, and Novikov's
  conjecture}.
\newblock {\em K-theory}, 4:303--337, 07 1991.

\bibitem{kasparov-skandalis03}
Gennadi Kasparov and Georges Skandalis.
\newblock {Groups Acting Properly on "Bolic" Spaces and the Novikov
  Conjecture}.
\newblock {\em Annals of Mathematics}, 158(1):165--206, 2003.

\bibitem{kasparov88}
G.G. Kasparov.
\newblock {Equivariant KK-theory and the Novikov conjecture}.
\newblock {\em Inventiones mathematicae}, 91(1):147--202, 1988.

\bibitem{lafforgue00}
V.~Lafforgue.
\newblock {A proof of property (RD) for cocompact lattices of SL(3,R) and SL(3,
  C)}.
\newblock {\em Journal of Lie Theory}, 10(2):255--267, 2000.

\bibitem{lafforgue02}
Vincent Lafforgue.
\newblock {K-th{\'e}orie bivariante pour les alg{\`e}bres de Banach et
  conjecture de Baum-Connes}.
\newblock {\em Inventiones mathematicae}, 149(1):1--95, Jul 2002.

\bibitem{mineyev-yu}
Igor Mineyev and Guoliang Yu.
\newblock {The Baum-Connes conjecture for hyperbolic groups}.
\newblock {\em Inventiones mathematicae}, 149(1):97--122, Jul 2002.

\bibitem{nica13}
Bogdan Nica.
\newblock {Proper isometric actions of hyperbolic groups on $L^p$-spaces}.
\newblock {\em Compositio Mathematica}, 149(5):773--792, 2013.

\bibitem{oberlin}
D.~M. Oberlin.
\newblock {$M_p(G)\neq M_q(G) (p^{-1}+q^{-1}=1)$}.
\newblock {\em Israel J. Math.}, 22(2):175--179, 1975.

\bibitem{phillips13}
N.~Christopher Phillips.
\newblock {Crossed products of Lp operator algebras and the K-theory of Cuntz
  algebras on Lp spaces}.
\newblock {\em https://arxiv.org/abs/1309.6406}, 2013.

\bibitem{pytlik}
Tadeusz Pytlik.
\newblock {A construction of convolution operators on free group}.
\newblock {\em Studia Math.}, (1):73--76, 1984.

\bibitem{sapir}
Mark Sapir.
\newblock {The rapid decay property and centroids in groups}.
\newblock {\em Journal of Topology and Analysis}, 07(03):513--541, 2015.

\bibitem{valette02}
Alain Valette.
\newblock {A Glimpse into Non-commutative Geometry: Property (RD). In:
  Introduction to the Baum-Connes Conjecture}.
\newblock {\em Lectures in Mathematics ETH Z{\"u}rich. Birkh{\"a}user, Basel},
  2002.

\bibitem{yu05}
Guoliang Yu.
\newblock {Hyperbolic groups admit proper affine isometric actions on lp
  -spaces}.
\newblock {\em Geometric {\&} Functional Analysis GAFA}, 15(5):1144--1151, Oct
  2005.

\end{thebibliography}

\Addresses

\end{document}